\documentclass[12pt,a4paper]{article}
\usepackage{amsfonts,amsmath,amssymb,amscd}
\usepackage{fancybox}
\usepackage[latin1]{inputenc}
\usepackage[usenames]{color}
\usepackage{graphicx}
\usepackage{float}
\usepackage{fullpage}
\usepackage{epsfig,amssymb}
\usepackage{enumerate}
\usepackage{subcaption}

\usepackage[normalem]{ulem}
\definecolor{avocado}{rgb}{0.34,0.51,0.01}

\newtheorem{thm}{Theorem}[section]
\newtheorem{cor}[thm]{Corollary}

\newtheorem{prop}[thm]{Proposition}
\newtheorem{defn}[thm]{Definition}
\newtheorem{rem}[thm]{Remark}

\begin{document}
	
\title{\textbf{Generalized Elastic Translating Solitons}}

\author{\'Alvaro P\'ampano}
\date{\today}

\maketitle 

\begin{abstract}
\noindent We study translating soliton solutions to the flow by powers of the curvature of curves in the plane. We characterize these solitons as critical curves for functionals depending on the curvature. More precisely, translating solitons to the flow by powers of the curvature are shown to be generalized elastic curves. In particular, focusing on the curve shortening flow, we deduce a new variational characterization of the grim reaper curve.
\\

\noindent{\emph{Keywords:} Curve Shortening Flow, Generalized Elastic Curves, Grim Reaper, Translating Solitons.}
\end{abstract}

\section{Introduction}

The study of geometric flows is a fruitful topic of investigation in the field of Differential Geometry. Mean curvature and inverse mean curvature flows have seen significant progress in the last decades, helping to bring insight to deep problems related to the theory of hypersurfaces (see, for instance, \cite{BN,BHW,CM,H1,HI,HI1,I}). In particular, translating solitons to these flows have been widely studied (see \cite{DLW,HIMW,KP} and references therein).

In the one dimensional case, considering the evolution of curves in the plane in which each point moves in a direction normal to the curve with speed involving a power of the curvature gives rise to some natural parabolic deformations \cite{A1,A2}. Beyond their mathematical interest, these curvature-driven flows also arise from applications in many different fields ranging from phase transition to image processing \cite{C,S,V}.

In this paper we consider the \emph{curvature-driven flow} of (possibly, convex) curves $X(\cdot, t):I\subseteq\mathbb{R}\longrightarrow\mathbb{R}^2$ given by
$$\frac{\partial X}{\partial t}(s,t)=\frac{1}{a}\left(\kappa^p(s,t)+b\right)N(s,t)\,.$$
Here $N(s,t)$ is the unit normal vector field along $X(s,t)$, $\kappa(s,t)$ is its curvature, and $a\neq 0$ and $b,p\in\mathbb{R}$ are fixed real constants. The case $p=1$ and $b=0$ is known as the \emph{curve shortening flow}, whose study was originated in the 1980s and now has an extensive literature. Some of the pioneering works are \cite{AL,G0,G1,GH,Gr} (see also \cite{A,Gi} and references therein).

A special class of solutions to the above partial differential equation is given by the self-similar solutions, which are those whose shape is independent of the evolution parameter $t$, possibly up to homotheties. For the curve shortening flow ($p=1$ and $b=0$) general self-similar solutions were classified by H. Halld\'orsson \cite{H}, while the classification of closed homothetic self-similar solutions was first achieved by U. Abresch and J. Langer \cite{AL}. More generally, homothetic self-similar solutions to the above described curvature-driven flow with $b=0$ have been completely classified by B. Andrews \cite{A1} and J. Urbas \cite{U}. 

The focus of this paper is on self-similar solutions that translate along a fixed direction only. These translational self-similar solutions are usually known as \emph{translating solitons}, which for the case $b=0$ were studied in \cite{NT}. In particular, for the curve shortening flow ($p=1$ and $b=0$) the only translating soliton is the \emph{grim reaper} curve (also known as the hairpin model). 

In the main result, we will show a characterization of the specific solutions to the above curvature-driven flow given by translating solitons in terms of critical curves for curvature dependent functionals.

\begin{thm}\label{intro} A smooth curve with nonconstant curvature $\gamma$ is a translating soliton to the curvature-driven flow
$$\frac{\partial X}{\partial t}(s,t)=\frac{1}{a}\left(\kappa^p(s,t)+b\right)N(s,t)\,$$
if and only if it is a critical curve, for compactly supported variations, for:
\begin{enumerate}[(i)]
	\item If $p\neq 1$,
	$$\mathbf{\Theta}_{p,\lambda}(\gamma):=\int_\gamma\left(\kappa^p+\lambda\right),$$
	where $\lambda=b\,(1-p)\in\mathbb{R}$.
	\item If $p=1$,
	$$\mathbf{\Theta}_{1,\lambda}(\gamma):=\int_\gamma\left(\kappa\log\kappa+\lambda\right),$$
	where $\lambda=-b\in\mathbb{R}$.
\end{enumerate}
\end{thm}

This result relates the topic of geometric flows with another central field of research in the area of Differential Geometry, namely, the analysis of functionals depending on the curvature.

The study of equilibria of functionals depending on the curvature goes back to the days of the Bernoulli family and L. Euler. Indeed, D. Bernoulli in a letter to L. Euler of 1738 proposed to investigate extrema of the functionals $\mathbf{\Theta}_{p\neq 1,\lambda}$, where the constant $\lambda\in\mathbb{R}$ may be understood as a Lagrange multiplier encoding the conservation of the length through the variation.

Clearly, the case $p=2$ recovers the pioneering bending energy employed in elasticity theory, whose critical curves were studied by J. Bernoulli in the 1690s and completely classified by L. Euler in an appendix to his book of 1744 (see \cite{T} for a historical background and references on this topic). The unconstrained critical curves, i.e., equilibria for $\mathbf{\Theta}_{2,0}$, also appear as the classical \emph{lintearia}, which represents the shape of a long cloth sheet full of water \cite{L} (extended notions of the \emph{lintearia} and their relations with equilibria for $\mathbf{\Theta}_{p,\lambda}$ can be found in \cite{LP}). Critical curves for $\mathbf{\Theta}_{2,\lambda}$ are usually referred as to elastic curves and so, the general case can be understood as an extension of this notion. Hence, critical curves for $\mathbf{\Theta}_{p,\lambda}$ are \emph{generalized elastic curves}.

Apart from the bending energy of curves, other functionals $\mathbf{\Theta}_{p,\lambda}$ for different choices of the energy parameters $p,\lambda\in\mathbb{R}$ have also been considered throughout the history. For instance, the case $p=1/2$ and $\lambda=0$ was considered by W. Blaschke in 1921 \cite{B}, who showed that critical curves with nonconstant curvature are catenaries by explicitly obtaining their curvatures in terms of their arc length parameters. W. Blaschke also considered the case $p=1/3$ and $\lambda=0$ in 1923 \cite{B}. This functional measures the equi-affine length for convex curves and its critical curves with nonconstant curvature are parabolas.

More recently, in \cite{LP0}, the critical curves with nonconstant curvature for the unconstrained cases $\mathbf{\Theta}_{p\neq 1,0}$ were related to a different variational problem that appears in the theory of weighted manifolds developed by M. Gromov \cite{G}, although it is not the usual weighted area functional that arises in the theory of mean curvature flows \cite{I}. To the contrary, these curves are the one dimensional analogue of (generalized) singular minimal surfaces \cite{D}. A particular case of these generalized singular minimal curves ($\alpha=-1/2$ and $\varpi=0$ in the functional $F_{\alpha,\varpi}$ of \cite{LP0}) gives rise to the \emph{brachistochrone} variational problem, which was posed by J. Bernoulli in 1696 and whose solution curves are cycloids (historical details can be found in \cite{BM}). It turns out that cycloids are also the critical curves for $\mathbf{\Theta}_{-1,0}$.

\begin{figure}[h!]
		\begin{subfigure}[b]{0.22\linewidth}
			\includegraphics[width=\linewidth]{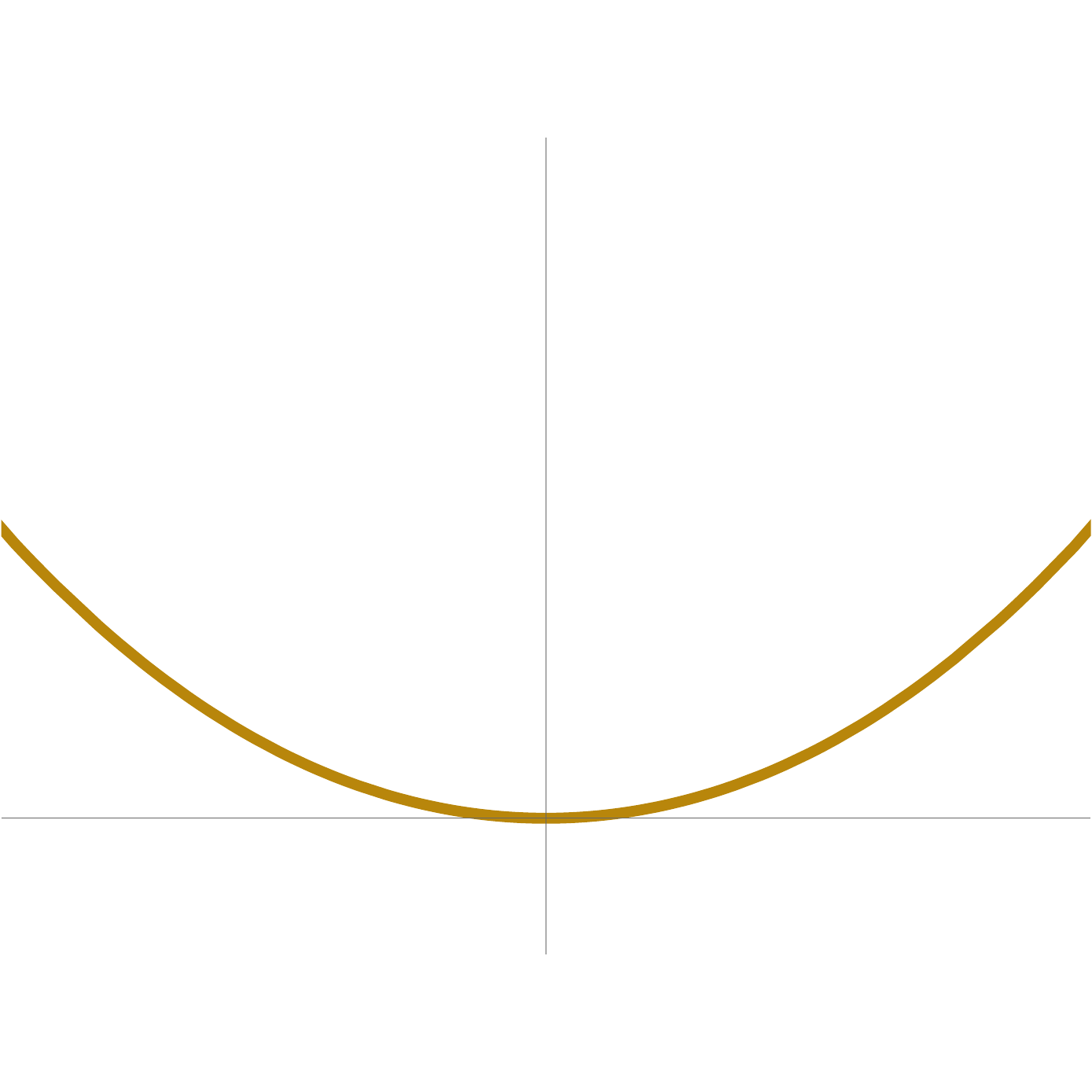}
			\caption{$\lambda=-0.5$}
		\end{subfigure}
		\quad
		\begin{subfigure}[b]{0.22\linewidth}
			\includegraphics[width=\linewidth]{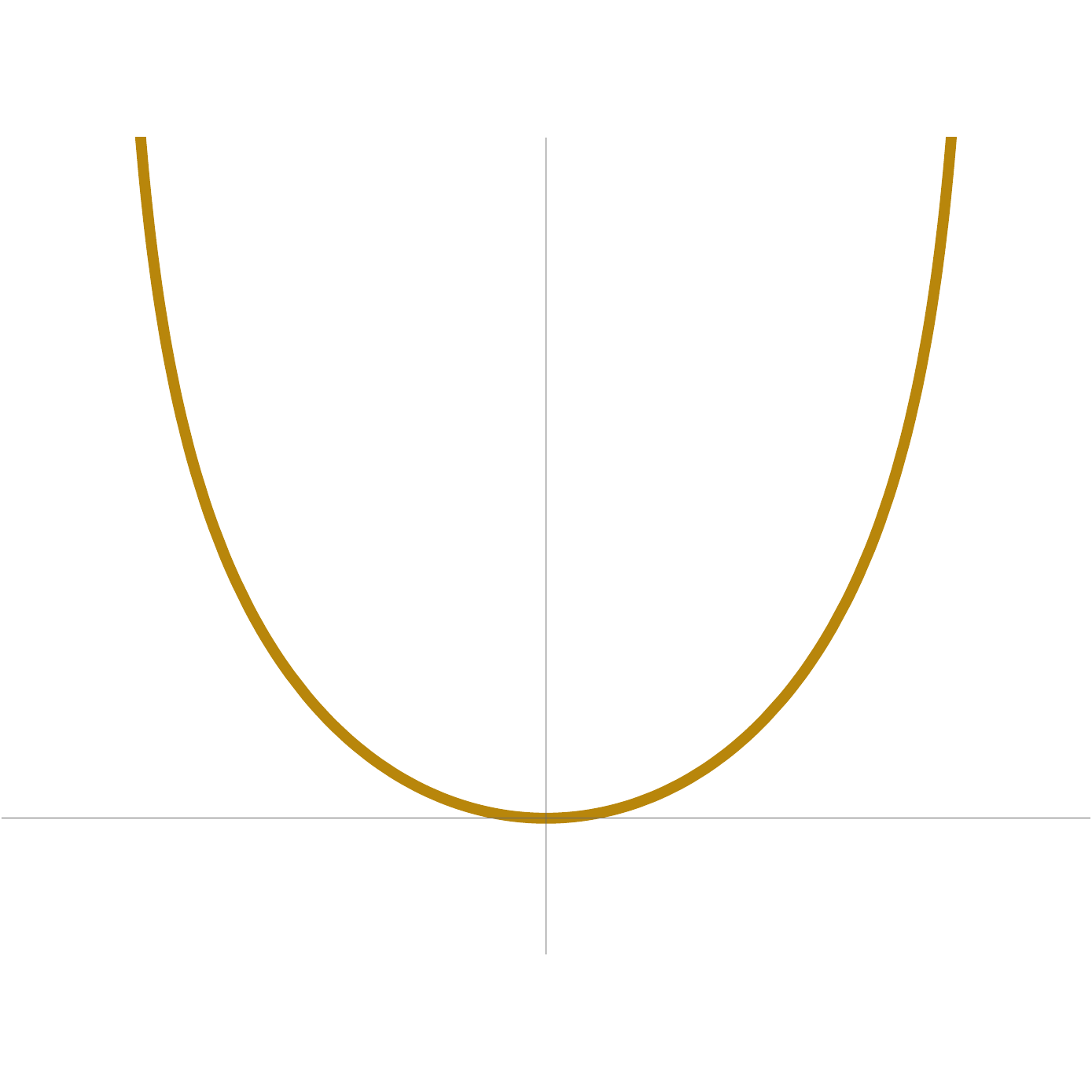}
			\caption{$\lambda=0$}
		\end{subfigure}
		\quad
		\begin{subfigure}[b]{0.22\linewidth}
			\includegraphics[width=\linewidth]{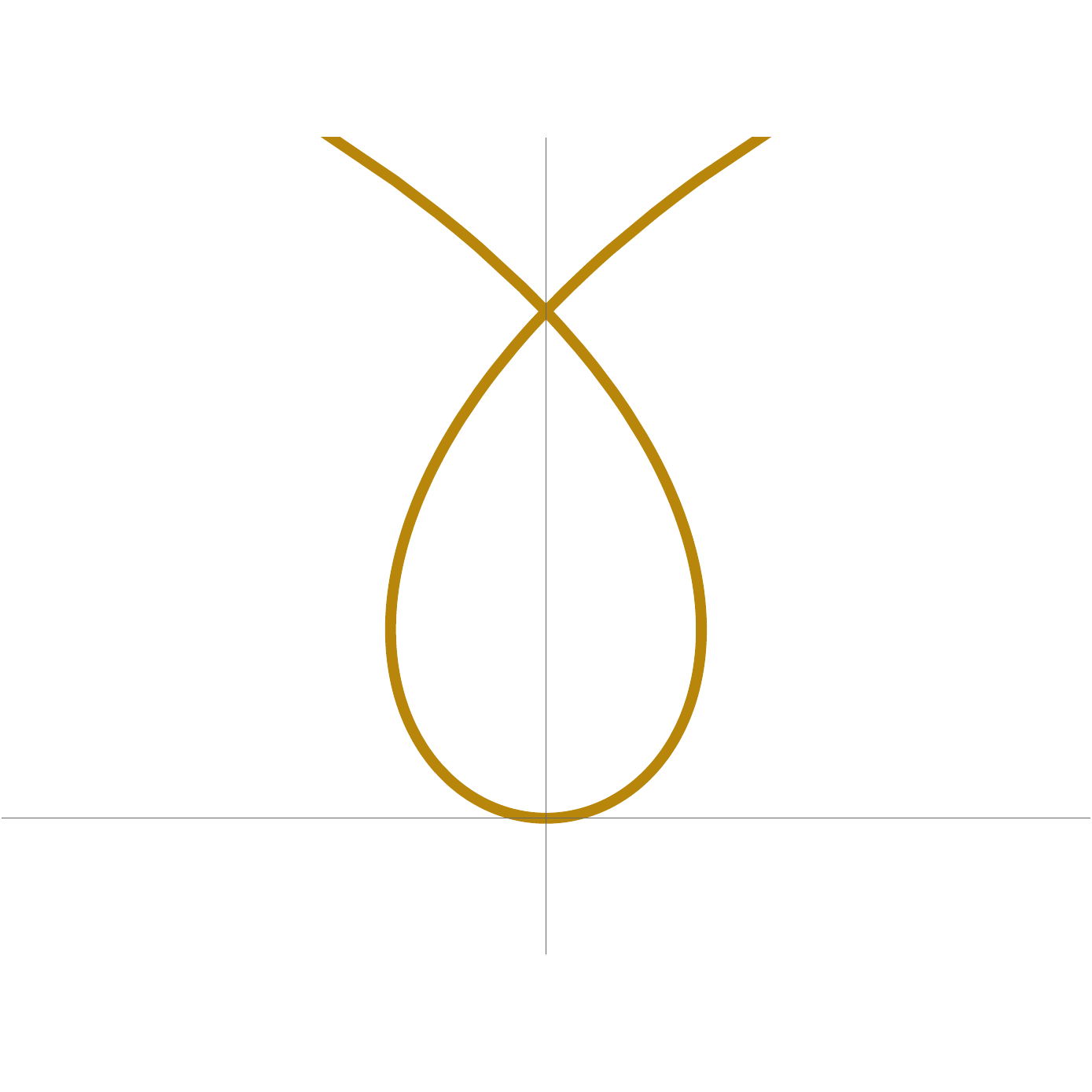}
			\caption{$\lambda=1$}
		\end{subfigure}
		\quad
		\begin{subfigure}[b]{0.22\linewidth}
			\includegraphics[width=\linewidth]{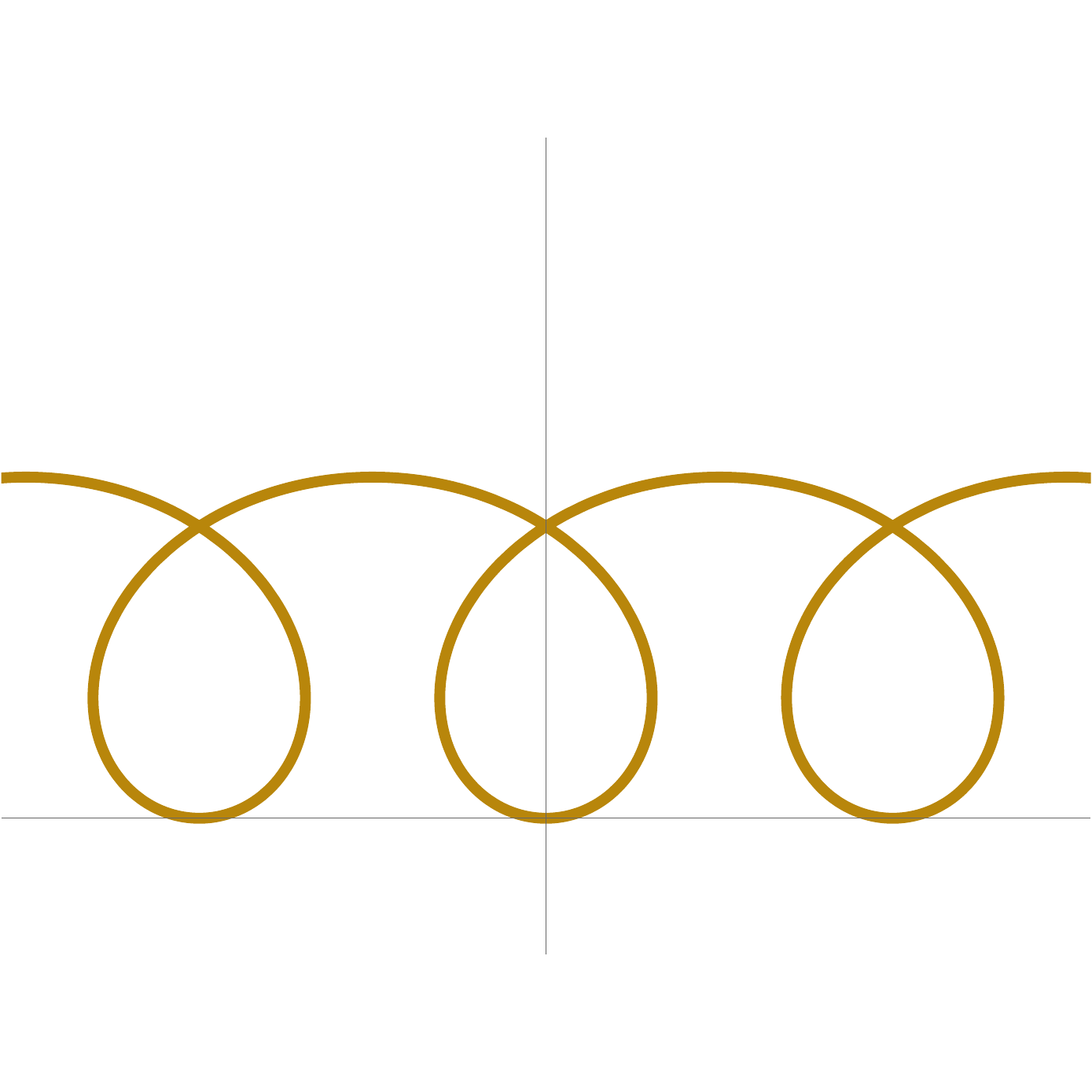}
			\caption{$\lambda=1.8$}
		\end{subfigure}
	\caption{Critical curves for the functional $\mathbf{\Theta}_{1,\lambda}$. These curves are the translating solitons to the curvature-driven flow given by $p=1$ and $b=-\lambda$, which evolve along the vertical axis. Figure (b) represents the grim reaper curve.}
	\label{Figure}
\end{figure}

In general, the critical curves with nonconstant curvature for $\mathbf{\Theta}_{p\neq 1,\lambda}$, were geometrically described in \cite{LP} and their general shape depicted in several figures (the case $p=1/2$ and $\lambda\in\mathbb{R}$ has also been recently considered in \cite{MP}). Observe that in \cite{LP} these curves are analyzed according to the value of $n=1/(p-1)$. To our best knowledge, the curvature depending energy functional $\mathbf{\Theta}_{1,\lambda}$ has not yet been studied, but the same analysis as for $\mathbf{\Theta}_{p\neq 1,\lambda}$ can be employed to geometrically describe the critical curves. For the sake of brevity this analysis will be omitted in this paper. Nonetheless, we show in Figure \ref{Figure} the shape of the critical curves with nonconstant curvature as the parameter $\lambda$ varies (c.f., \cite{BO}).

The relation introduced in Theorem \ref{intro} allows us to employ the tools and techniques of the geometric variational problems associated with $\mathbf{\Theta}_{p,\lambda}$ to the study of translating solitons to the above introduced curvature-driven flows. For instance, it is known that critical curves for $\mathbf{\Theta}_{p,\lambda}$ can be parameterized locally in terms of just one quadrature (see Proposition \ref{coord} and Remark \ref{rem} below). In the specific case of $\mathbf{\Theta}_{1,\lambda}$, since the curvature of the critical curves can be explicitly obtained from the corresponding Euler-Lagrange equation, this parameterization can indeed be given in terms of elementary functions. This coincides with \cite{BO}, although the approach to the problem is different.

\section{Flow by Powers of the Curvature}

Let $X:I\subseteq\mathbb{R}\longrightarrow \mathbb{R}^2$ be the smooth immersion of a curve $\gamma$ parameterized by the arc length parameter $s\in I$. Denote by $T(s):=X'(s)$ the unit tangent vector field along the curve $\gamma$, where $\left(\,\right)'$ denotes the derivative with respect to the arc length parameter $s$, and define the unit normal vector field $N(s)$ along $X(s)$ to be the counter-clockwise rotation of $T(s)$ through an angle $\pi/2$. In this setting, the (signed) \emph{curvature} $\kappa(s)$ of $\gamma$ is defined by the Frenet-Serret equation
$$T'(s)=\kappa(s)N(s)\,.$$
A curve $\gamma$ is a \emph{geodesic} (i.e., a straight line) if its curvature $\kappa(s)$ is identically zero, and it is \emph{convex} if $\kappa(s)>0$ everywhere.

Let $\gamma_t$, $t\in [0,t_o)$, be a family of curves given by smooth immersions $X_t\equiv X(\cdot,t):I\subseteq\mathbb{R}\longrightarrow \mathbb{R}^2$ where $I$ is some open interval in $\mathbb{R}$. The curves $\gamma_t$ are said to evolve under a \emph{curvature-driven flow} if
\begin{equation}\label{flow}
	\frac{\partial X}{\partial t}(s,t)=\frac{1}{a}\left(\kappa^p(s,t)+b\right)N(s,t)\,,
\end{equation}
for all $(s,t)\in I\times [0,t_o)$, where $s\in I$ is the arc length parameter of the curve $\gamma_t$, $\kappa(s,t)$ is its curvature, $N(s,t)$ is the unit normal vector field along $X_t(s)$, and $a\neq 0$ and $b,p\in\mathbb{R}$ are fixed real constants. If $p\in\mathbb{R}\setminus\mathbb{N}$ (we consider $0\in\mathbb{N}$), any smooth solution of \eqref{flow} must have $\kappa>0$ and so, in that case, possibly after a change of orientation, we will restrict ourselves to convex curves.

In this paper, we focus on a special type of flow by imposing that the curves $\gamma_t$ translate along a direction. The general form of solutions to \eqref{flow} translating along a direction is
\begin{equation}\label{motion}
	X(s,t)=X(s)+t V\,,
\end{equation}
for a constant vector $V\in\mathbb{R}^2$. Consequently, from \eqref{flow}, they satisfy the equation
\begin{equation}\label{curvt}
	\kappa^p(s,t)+b=a\,\langle N(s,t),V\rangle\,.
\end{equation}
Indeed, up to tangential diffeomorphisms, equation \eqref{curvt} is also a sufficient condition for \eqref{motion} to be a solution of \eqref{flow}. In particular, since the curves $\gamma_t$ translate along a direction, the curvature $\kappa(s,t)$ and the unit normal $N(s,t)$ are independent of the parameter $t\in[0,t_o)$ and so they are the same as those of the initial condition $X(s,0)=X(s)$. For simplicity, we then introduce the following definition.

\begin{defn} An arc length parameterized smooth curve $X(s)\subset\mathbb{R}^2$ is a translating soliton to the curvature-driven flow \eqref{flow} if
	\begin{equation}\label{curv}
		\kappa^p(s)+b=a\,\langle N(s),V\rangle\,,
	\end{equation}
where $\kappa(s)$ is the curvature of $X(s)$, $N(s)$ is its unit normal and $V\in\mathbb{R}^2$ is a constant vector.
\end{defn}

We will consider now that a translating soliton has constant curvature $\kappa(s)=\kappa_o$. In this case, the arc length parameterized curve $X(s)$ may be a part of a straight line, if $\kappa_o=0$, or a part of a circle, if $\kappa_o\neq 0$. In the latter case, the left hand side of \eqref{curv} is constant, but the right hand side is not. Thus, we obtain that this case is not possible. In other words, the only translating solitons to the flow \eqref{flow} with constant curvature, whenever they exist (i.e., for suitable choices of the parameters $b,p\in\mathbb{R}$), are straight lines evolving in the constant direction given by their tangent vector.

In view of the above observations, to obtain non-trivial solutions to the curvature-driven flow \eqref{flow} translating along a direction, it is enough to solve the prescribed curvature equation \eqref{curv} for non-geodesics.

\section{Curvature Energy Functionals} 

Before stating and proving the main result of the paper, we briefly recall in this section the general theory for curvature energies. Consider a functional 
$$\mathbf{\Theta}(\gamma):=\int_\gamma P\left(\kappa\right),$$
where $P(\kappa)$ is a smooth function defined on an adequate domain, $\kappa$ is the curvature of the curve $\gamma$ and $\mathbf{\Theta}$ is acting on the space of curves immersed in $\mathbb{R}^2$. Regardless of the possible boundary conditions, a critical curve for the functional $\mathbf{\Theta}$ must satisfy the associated Euler-Lagrange equation
\begin{equation}\label{EL}
	\frac{d^2 \dot{P}}{ds^2}+\kappa^2\dot{P}-\kappa P=0\,,
\end{equation}
where $\dot{P}$ is the derivative of $P$ with respect to $\kappa$ and $s$ represents the arc length parameter of the curve. Equation \eqref{EL} is also the sufficient condition to characterize critical curves for $\mathbf{\Theta}$ for compactly supported variations. If $d\dot{P}/ds\neq 0$, the Euler-Lagrange equation \eqref{EL} can be integrated, obtaining
\begin{equation}\label{fi}
	\left(\frac{d \dot{P}}{ds}\right)^2+\left(\kappa\dot{P}-P\right)^2=d\,,
\end{equation}
for some positive constant of integration $d\in\mathbb{R}^+$. Observe that the condition $d\dot{P}/ds=\ddot{P}\kappa'\neq 0$ implies that the curvature $\kappa(s)$ of the critical curve for $\mathbf{\Theta}$ cannot be constant and that $P(\kappa)$ cannot be an affine function of the curvature.

For our purposes, it is convenient to use a special parameterization of critical curves for $\mathbf{\Theta}$ with nonconstant curvature. Although the following result is well-known in the theory of curvature dependent functionals, we include a proof for the sake of completeness.

\begin{prop}\label{coord}
	Let $d>0$ and $X(s)\subset\mathbb{R}^2$ be an arc length parameterized smooth curve with nonconstant curvature $\kappa(s)$. Then, $\kappa(s)$ is a solution of \eqref{fi} if and only if there exists a coordinate system in $\mathbb{R}^2$ such that
	\begin{equation}\label{param}
		X(s)=\frac{1}{\sqrt{d\,}}\left(\int\left(\kappa\dot{P}\left(\kappa\right)-P\left(\kappa\right)\right)ds,\dot{P}\left(\kappa\right)\right),
	\end{equation}
	where $\kappa\equiv\kappa(s)$.
\end{prop}
\textit{Proof.} The forward implication follows from standard computations involving Killing vector fields along curves (see \cite{LS} for a definition of these fields). Let $X(s)=\left(x_1(s),x_2(s)\right)$ be the smooth immersion of a curve $\gamma$ whose nonconstant curvature $\kappa\equiv\kappa(s)$ satisfies \eqref{fi}. Then, the vector field
\begin{equation}\label{J}
	\mathcal{J}(s)=\left(\kappa\dot{P}\left(\kappa\right)-P\left(\kappa\right)\right)T(s)+\frac{d\dot{P}}{ds}\left(\kappa\right)N(s)\,,
\end{equation}
is a Killing vector field along the curve $\gamma$. Killing vector fields along curves can be uniquely extended to Killing vector fields on $\mathbb{R}^2$. For convenience, we denote the extension of \eqref{J} also by $\mathcal{J}$. Observe that equation \eqref{fi} can be simply rewritten as $\langle \mathcal{J},\mathcal{J}\rangle=d$ and so the extension of $\mathcal{J}$ to $\mathbb{R}^2$ is a translational Killing vector field. After a rigid motion if necessary, we will assume $\mathcal{J}=\sqrt{d}\,\partial_{x_1}$. It then follows that
\begin{equation}\label{x2}
	\kappa\dot{P}\left(\kappa\right)-P\left(\kappa\right)=\langle T,\mathcal{J}\rangle=\sqrt{d}\,x_1'(s)\,,
\end{equation}
and so, from the arc length parameterization of the curve $\gamma$ given by the immersion $X(s)$ and from \eqref{fi}, we conclude that
\begin{equation}\label{x1}
	x_2'(s)=\frac{1}{\sqrt{d}}\,\frac{d\dot{P}}{ds}\left(\kappa\right).
\end{equation}
After integrating \eqref{x2} and \eqref{x1} we deduce that, up to rigid motions, the arc length parameterized curve $X(s)=\left(x_1(s),x_2(s)\right)$ whose curvature $\kappa\equiv\kappa(s)$ satisfies \eqref{fi} is given by \eqref{param}. This proves the forward implication.

We now show the converse. Assume that the arc length parameterization $X(s)$ of a curve $\gamma$ is given by \eqref{param} for some $d>0$ and where $\kappa\equiv\kappa(s)$ denotes its curvature. Differentiating with respect to the arc length parameter we obtain the unit tangent vector $T(s)=X'(s)$. The condition that this vector is unitary shows that \eqref{fi} holds, which proves the statement. \hfill$\square$

\begin{rem}\label{rem} If the first integral \eqref{fi} of the Euler-Lagrange equation is employed to make a change of variable in the integral of the first component of the parameterization \eqref{param}, we locally obtain a new parameterization (with the curvature as the parameter) depending on just one quadrature.
\end{rem}

The parameterization of critical curves for $\mathbf{\Theta}$ obtained in above result will play an essential role in the proof of the characterization of translating solitons to the curvature-driven flow \eqref{flow}.

\section{Variational Characterization}

In this section we will prove the main result of the paper, in which we show that translating solitons to the curvature-driven flow \eqref{flow} are characterized as critical curves with respect to compactly supported variations for suitable energy functionals depending on the curvature. The following theorem is equivalent to Theorem \ref{intro} of the Introduction.

\begin{thm}\label{main}
	Let $X(s)\subset\mathbb{R}^2$ be an arc length parameterized smooth curve with nonconstant curvature. Then $X(s)$ is a translating soliton to the curvature-driven flow \eqref{flow} if and only if its curvature $\kappa\equiv\kappa(s)$ is a solution to the Euler-Lagrange equation associated with:
	\begin{enumerate}[(i)]
		\item If $p\neq 1$,
		\begin{equation}\label{energy}
			\mathbf{\Theta}_{p,\lambda}(\gamma):=\int_\gamma\left(\kappa^p+\lambda\right),
		\end{equation}
		where $\lambda=b\,(1-p)$ is a real constant.
		\item If $p=1$,
		\begin{equation}\label{energy1}
			\mathbf{\Theta}_{1,\lambda}(\gamma):=\int_\gamma\left(\kappa\log\kappa+\lambda\right),
		\end{equation}
		where $\lambda=-b$ is a real constant.
	\end{enumerate}
\end{thm}
\textit{Proof.} Assume first that $X(s)$ is the arc length parameterization of a translating soliton to the flow \eqref{flow} with nonconstant curvature. By definition, this means that the nonconstant curvature $\kappa(s)$ of $X(s)$ satisfies \eqref{curv}. After a suitable rigid motion, we can assume without loss of generality that the vector $V$ is a multiple of $(0,1)$. Moreover, possibly varying the value of the parameter $a\neq 0$ (which will not appear in our characterization) we may assume $V=(0,1)$. It then follows from \eqref{curv} that
$$\langle N(s), (0,1)\rangle=\frac{1}{a}\left(\kappa^p+b\right),$$
and so, from the definition of the unit normal $N(s)$ as a counter-clockwise rotation of the tangent $T(s)$,
$$\langle T(s), (1,0)\rangle=\frac{1}{a}\left(\kappa^p+b\right).$$
If the translating soliton $X(s)$ were to be parameterized as \eqref{param}, then the following relation should hold
\begin{equation}\label{ODE}
	\frac{1}{a}\left(\kappa^p+b\right)=\frac{1}{\sqrt{d}\,}\left(\kappa \dot{P}(\kappa)-P(\kappa)\right).
\end{equation}
This represents an equation independent of the arc length parameter $s$ and so it can be understood as an ordinary differential equation in $P(\kappa)$. The ordinary differential equation \eqref{ODE} can be explicitly solved obtaining
$$P(\kappa)=\frac{\sqrt{d}}{a(p-1)}\,\kappa^p+\mu\,\kappa-\frac{b\,\sqrt{d}}{a}\,,$$
if $p\neq 1$, or
$$P(\kappa)=\frac{\sqrt{d}}{a}\,\kappa \log\kappa+\mu\,\kappa-\frac{b\sqrt{d}}{a}\,,$$
if $p=1$. In both cases $\mu\in\mathbb{R}$ is a constant of integration. However, it follows from the general equation \eqref{EL} that the linear term in the above expressions does not affect the Euler-Lagrange equation and so, we may take without loss of generality that $\mu=0$. Furthermore, any multiple of $P(\kappa)$ will clearly give the same Euler-Lagrange equation \eqref{EL} and so, we can take the functions $P(\kappa)$ to be those of the statement for the corresponding values of $\lambda$.

It remains to prove that for those choices of $P(\kappa)$ the translating soliton $X(s)$ and the arc length parameterized curve $\widetilde{X}(s)$ whose curvature satisfies the Euler-Lagrange equation \eqref{EL} locally coincide. Observe that, from Proposition \ref{coord}, the curve $\widetilde{X}(s)$ can be parameterized as \eqref{param} (replacing $\kappa$ by $\widetilde{\kappa}$, the curvature of $\widetilde{X}(s)$). Then, perhaps up to orientation, the curvature $\widetilde{\kappa}(s)$ of $\widetilde{X}(s)$ locally coincides with that of the translating soliton $X(s)$. The Fundamental Theorem of Curves then implies that locally $X(s)=\widetilde{X}(s)$ and so the curvature of the translating soliton also satisfies \eqref{EL} for the corresponding value of $P(\kappa)$. This finishes the proof of the forward implication, namely, translating solitons to \eqref{flow} satisfy the Euler-Lagrange equation associated with either \eqref{energy} or \eqref{energy1}, respectively.

The converse follows directly. Indeed, if $X(s)$ is an arc length parameterized curve whose curvature satisfies the Euler-Lagrange equation associated with either $\mathbf{\Theta}_{p,\lambda}$, then from Proposition \ref{coord} it can be parameterized as \eqref{param} for the corresponding value of $P(\kappa)$. Differentiating with respect to the arc length parameter we obtain the unit tangent $T(s)=X'(s)$ and after the suitable counter-clockwise rotation through an angle $\pi/2$, we deduce the expression for the unit normal $N(s)$. It is then a simple verification to check that \eqref{curv} holds, and so $X(s)$ is a translating soliton to the curvature-driven flow \eqref{flow}. \hfill$\square$
\\

The above result yields several unexpected characterizations of translating solitons. In fact, combining the known results on the field of functionals depending on the curvature with those of translating solitons, we obtain the following particular cases:
\begin{enumerate}[(i)]
	\item Case $p=-1$ and $b=0$. For these values, equation \eqref{flow} represents the inverse mean curvature flow, whose translating solitons were proven to be cycloids in \cite{DLW}. Cycloids were classically known to solve the \emph{brachistochrone} problem. In \cite{LP0}, it was shown that cycloids are also equilibria for the total radius of curvature, i.e., critical curves for $\mathbf{\Theta}_{-1,0}$, which coincides with Theorem \ref{main}. The fact that cycloids are solutions to the variational problem associated with $\mathbf{\Theta}_{-1,0}$ was likely already known.
	\item Case $p=0$. Translating solitons to the flow \eqref{flow} with $p=0$ are necessarily straight lines. Indeed, equation \eqref{curv} shows that the angle between the unit normal $N(s)$ and a vector $V\in\mathbb{R}^2$ must be constant, and so the angle between $V$ and the unit tangent $T(s)$ must also be constant. This agrees with the result of Theorem \ref{main} since, for $p=0$, the functional $\mathbf{\Theta}_{0,\lambda}$ is the length functional whose only critical curves are straight lines.
	\item Case $p=1/3$ and $b=0$. W. Blaschke proved that critical curves with nonconstant curvature for $\mathbf{\Theta}_{1/3,0}$ are parabolas \cite{B} and so, it follows from Theorem \ref{main} that the only translating solitons to the flow \eqref{flow} with $p=1/3$ and $b=0$ are parabolas.
	\item Case $p=1/2$ and $b=0$. Translating solitons to \eqref{flow} with these values are catenaries since, according to Theorem \ref{main}, they must be equilibria for $\mathbf{\Theta}_{1/2,0}$. W. Blaschke obtained that the only critical curves, other than straight lines, are catenaries \cite{B}.
	\item Case $p=2$. The classical elastic curves of Euler-Bernoulli are the translating solitons to \eqref{flow} for $p=2$.
\end{enumerate}

In particular, special mention is deserved by the case $p=1$ for which Theorem \ref{main} gives a new characterization of grim reaper curves in terms of critical curves with respect to compactly supported variations for an energy functional depending on the curvature. To see this, consider the curve shortening flow ($p=1$ and $b=0$ in \eqref{flow}) given by
$$\frac{\partial X}{\partial t}(s,t)=\frac{1}{a}\kappa(s,t)N(s,t)\,.$$
It is well-known that the non-trivial translating solitons are the grim reaper curves. Combining this fact with the result of Theorem \ref{main} we conclude with a new variational characterization of these curves.

\begin{cor} An arc length parameterized nongeodesic curve $X(s)\subset\mathbb{R}^2$ is a grim reaper if and only if its curvature $\kappa\equiv\kappa(s)$ is a solution to the Euler-Lagrange equation associated with 
$$\mathbf{\Theta}_{1,0}(\gamma):=\int_\gamma \kappa\log\kappa\,.$$
\end{cor}

\section*{Appendix. Logarithmic Curvature Flow}

In this appendix, we apply the same technique to a different flow, namely, the \emph{logarithmic curvature flow} given by (c.f., \cite{CW})
\begin{equation}\label{log}
	\frac{\partial X}{\partial t}(s,t)=\frac{1}{a}\left(\log\kappa(s,t)+b\right)N(s,t)\,,
\end{equation}
where $\kappa(s,t)$ is the curvature of the curve $\gamma_t$, $N(s,t)$ is its unit normal vector field and $a\neq 0$ and $b\in\mathbb{R}$ are fixed constants. Of course, the curves under consideration are convex, i.e., $\kappa(s,t)>0$ holds everywhere along the curve.

In this setting, and after the same observations of Section 2, we will say that an arc length parameterized smooth curve $X(s)\subset\mathbb{R}^2$ is a translating soliton to the logarithmic curvature flow \eqref{log} if 
\begin{equation}\label{curvlog}
	\log(\kappa(s))+b=a\,\langle N(s),V\rangle\,,
\end{equation}
for a constant vector $V\in\mathbb{R}^2$. We then obtain the following variational characterization.

\begin{thm}
	Let $X(s)\subset\mathbb{R}^2$ be an arc length parameterized smooth curve. Then $X(s)$ is a translating soliton to the logarithmic curvature flow \eqref{log} if and only if its curvature $\kappa\equiv\kappa(s)$ is a nonconstant solution of the Euler-Lagrange equation associated with
	\begin{equation}\label{energylog}
		\widetilde{\mathbf{\Theta}}_\lambda(\gamma):=\int_\gamma\left(\log\kappa+\lambda\right),
	\end{equation}
	where $\lambda=b+1$ is a real constant.
\end{thm}
\textit{Proof.} The proof follows the same steps as that of Theorem \ref{main}. Hence, we will omit it. \hfill$\square$
\\

The curvature dependent functional $\widetilde{\mathbf{\Theta}}_\lambda$ was studied in \cite{LP}, where the critical curves were described. Their shape is similar to those of Figure 3 of \cite{LP}.

\begin{flushleft}
\'Alvaro P{\footnotesize \'AMPANO}\\
Department of Mathematics and Statistics, Texas Tech University, Lubbock, TX, 79409, USA\\
E-mail: alvaro.pampano@ttu.edu
\end{flushleft}

\end{document}